\documentclass[11pt,a4paper]{article}

\usepackage{inputenc}
\usepackage{amsmath}
\usepackage{bm}
\usepackage{bbold}
\usepackage{amsthm}

\usepackage{curves}

\usepackage{hyperref}

\title{Growth Rate of the State Vector in a Generalized Linear Stochastic System with Symmetric Matrix\thanks{Journal of Mathematical Sciences, 2007, Vol.~147, N~4, pp.~6924--6928.}
} 

\author{N.~K.~Krivulin\thanks{Faculty of Mathematics and Mechanics, St.~Petersburg State University, 28 Universitetsky Ave., St.~Petersburg, 198504, Russia, nkk@math.spbu.ru.}\thanks{The work was partially supported by the Russian Foundation for Basic Research (grant no. 06-01-00763).}}

\date{}

\begin{document}

\maketitle

\begin{abstract}
The mean growth rate of the state vector is evaluated for a generalized linear stochastic second-order system with a symmetric matrix. Diagonal entries of the matrix are assumed to be independent and exponentially distributed with different means, while the off-diagonal entries are equal to zero.
%
\end{abstract}

\section{Introduction}

The analysis of actual systems in engineering, economics, manufacturing, and other areas often involves generalized linear dynamical models which take the form
$$
\bm{z}(k)=A(k)\otimes\bm{z}(k-1),
$$
where $ A(k) $ is a random state transition matrix, $ \bm{z}(k) $ is a state vector, and $ \otimes $ denotes matrix-vector multiplication in some idempotent algebra \cite{Baccelli1993Synchronization,Maslov1994Idempotent}.

In the analysis, one is normally interested in evaluation of the mean (asymptotic) growth rate of the system state vector
$$
\lambda
=
\lim_{k\to\infty}\frac{1}{k}\|\bm{z}(k)\|,
$$
where $ \|\cdot\| $ is an idempotent counterpart to the usual vector norm.

In many cases, it is not difficult to prove existence for the above limit (e.g., with the ergodic theorem presented in \cite{Kingman1973Subadditive}). At the same time, evaluation of the limit may appear to be a difficult problem even for quite simple systems. The related results include the solution obtained in \cite{Olsder1990Discrete} for the second-order system with the matrix having all its entries independent and exponentially distributed with unit mean. In \cite{Jeanmarie1994Analytical}, another second-order system is examined which has a symmetric matrix such that its diagonal entries are independent and exponentially distributed with unit mean, whereas the off-diagonal entries are equal to zero.

For systems of arbitrary order $ n $, there are results known for the case when the entries of the matrix are all independent and identically distributed with either normal \cite{Cohe1998Subadditivity} or discrete uniform \cite{Olsder1990Discrete} probability distributions. In \cite{Krivulin2005Growth}, a general solution is presented for a system with a triangular matrix provided that its random entries have arbitrary distributions with finite mean and variance, and may be dependent.

The purpose of this paper is to evaluate the mean growth rate of the state vector for a second-order system which somewhat generalizes the system with symmetric matrix examined in \cite{Jeanmarie1994Analytical}. In the paper, we assume that the off-diagonal entries of the matrix are still equal to zero, whereas the diagonal entries are independent, and have exponential distributions with parameters $ \mu $ and $ \nu $, $ \mu\ne\nu $.

\section{Stochastic Linear Dynamical System}

Consider a dynamical system such that its evolution is described by the equation
\begin{equation}\label{E-LSDS}
\bm{z}(k)=A(k)\otimes\bm{z}(k-1),
\end{equation}
where
$$
A(k)
=
\begin{pmatrix}
	\alpha_{k}	& 0 \\
	0						& \beta_{k}
\end{pmatrix},
\qquad
\bm{z}(k)
=
\begin{pmatrix}
	x(k) \\
	y(k)
\end{pmatrix},
\qquad
\bm{z}(0)
=
\begin{pmatrix}
	0 \\
	0
\end{pmatrix},
$$
and symbol $ \otimes $ stands for matrix-vector multiplication in idempotent algebra \cite{Baccelli1993Synchronization,Maslov1994Idempotent} with scalar operations $ \oplus $ and $ \otimes $ defined respectively as maximum and addition.

Equation \eqref{E-LSDS} can be represented as the system of scalar equations
\begin{align*}
x(k)&=\alpha_{k}\otimes x(k-1)\oplus y(k-1), \\
y(k)&= x(k-1)\oplus\beta_{k}\otimes y(k-1),
\end{align*}
or, with ordinary notation, as 
\begin{align*}
x(k)&=\max(\alpha_{k}+x(k-1),y(k-1)), \\
y(k)&=\max(x(k-1),\beta_{k}+y(k-1)).
\end{align*}

Suppose the sequences $ \{\alpha_{k} \} $ and $ \{\beta_{k}\} $ each involve independent and identically distributed random variables, and $ \alpha_{k} $ and $ \beta_{l} $ are independent for any $ k,l $.

Further assume that $ \alpha_{k} $ and $ \beta_{k} $ have the exponential probability distribution with respective parameters $ \mu $ and $ \nu $. Denote their related distribution functions by
$$
F(t)
=
\max(0,1-e^{-\mu t}),
\qquad
G(t)
=
\max(0,1-e^{-\nu t}),
$$
and their density functions by $ f(t) $ and $ g(t) $.

Note that if $ \mu=\nu $, then the system can easily be reduced to that with $ \mu=\nu=1 $ examined in \cite{Jeanmarie1994Analytical}.

Among the performance measures of interest for system \eqref{E-LSDS} is the mean growth rate of system state vector
\begin{equation}
\lambda=\lim_{k\to\infty}\frac{1}{k}\|\bm{z}(k)\|,\label{E-lambda}
\end{equation}
where the symbol of norm is taken in an idempotent algebra sense,
$$
\|\bm{z}(k)\|=x(k)\oplus y(k)=\max(x(k),y(k)).
$$

\section{Mean Growth Rate of State Vector}

By applying the ergodic theorem in \cite{Kingman1973Subadditive}, it is not difficult to show that if the entries of the matrix $ A(k) $ are nonnegative and have finite mean values, then limit \eqref{E-lambda} exists with probability one. Moreover, there also exists the limit
$$
\lim_{k\to\infty}\frac{1}{k}\mathbb{E}\|\bm{z}(k)\|=\lambda.
$$

In order to evaluate the limit, first introduce the random variables
$$
Z(k)
=
\|\bm{z}(k)\|
-
\|\bm{z}(k-1)\|,
\qquad
Y(k)
=
y(k)-x(k),
$$
and note that $ \|\bm{z}(k)\|=Z(1)+\cdots+Z(k) $ and $ Y(0)=0 $.

Consider the distribution functions
$$
\Phi_{k}(t)
=
\mathbb{P}\{Z(k)<t\},
\qquad
\Psi_{k}(t)
=
\mathbb{P}\{Y(k)<t\}.
$$

The function $ \Phi_{k} $ can be represented in the form
\begin{multline*}
\Phi_{k}(t)
=
\mathbb{P}\{\max(x(k-1)+\alpha_{k},y(k-1)+\beta_{k})-\max(x(k-1),y(k-1))<t\} =\\
=
\mathbb{P}\{\max(\alpha_{k},Y(k-1)+\beta_{k})-\max(0,Y(k-1))<t\},
\end{multline*}
which leads to the equation
$$
\Phi_{k}(t)
=
G(t)\left(1-\int_{0}^{\infty}\Psi_{k-1}(u)f(u+t)du\right)
+
F(t)\int_{0}^{\infty}\Psi_{k-1}(-u)g(u+t)du.
$$

Suppose the sequence of functions $ \Psi_{k} $ tends to a distribution function $ \Psi $ as $ k\to\infty $. (Below we show that this assumption is valid for the system under study.)

Application of Lebesgue's dominated convergence theorem yields the conclusion that the sequence $ \Phi_{k} $ converges to the distribution function
\begin{equation}
\Phi(t)
=
G(t)\left(1-\int_{0}^{\infty}\Psi(u)f(u+t)du\right)
+
F(t)\int_{0}^{\infty}\Psi(-u)g(u+t)du\label{E-Phi}
\end{equation}
of some random variable $ Z $, and that $ Z(k)\to Z $ in distribution.

Moreover, $ \mathbb{E}[Z(k)] $ converges to $ \mathbb{E}[Z] $, and consequently,
\begin{equation}
\lambda
=
\lim_{k\to\infty}\frac{1}{k}\mathbb{E}\|\bm{z}(k)\|
=
\lim_{k\to\infty}\frac{1}{k}\sum_{m=1}^{k}\mathbb{E}[Z(m)]
=
\mathbb{E}[Z].\label{E-lambda1}
\end{equation}

\section{A Recursive Equation for $ \Psi_{k} $}

With the law of total probability, one can represent $ \Psi_{k} $ as
$$
\Psi_{k}(t)
=
\int_{0}^{\infty}\int_{0}^{\infty}
\mathbb{P}\{Y(k)<t|\alpha_{k}=u,\beta_{k}=v\}f(u)g(v)dudv.
$$

Considering that
\begin{multline*}
\mathbb{P}\{Y(k)<t|\alpha_{k}=u,\beta_{k}=v\} =\\
=
\mathbb{P}\{\max(x(k-1),y(k-1)+v)-\max(x(k-1)+u,y(k-1))<t\} =\\
=
\mathbb{P}\{\max(0,Y(k-1)+v)-\max(u,Y(k-1))<t\} =\\
=
\begin{cases}
1-\Psi_{k-1}(-t),	& \text{if $ u\leq-t $, $ v<t $}, \\
1,	& \text{if $ u>-t $, $ v<t $}, \\
0,	& \text{if $ u\leq-t $, $ v\geq t $}, \\
\Psi_{k-1}(u-v+t),	& \text{if $ u>-t $, $ v\geq t $},
\end{cases}
\end{multline*}
we arrive at the recursive equation
$$
\Psi_{k}(t)
=
\begin{cases}
\int\limits_{0}^{\infty}\int\limits_{0}^{\infty}\Psi_{k-1}(u-v)f(u-t)g(v)dudv,			& \text{if $ t\leq0 $}, \\
\int\limits_{0}^{\infty}\int\limits_{0}^{\infty}\Psi_{k-1}(u-v)f(u)g(v+t)dudv+G(t),	& \text{if $ t>0 $}.
\end{cases}
$$

Taking into account exponential distribution of the diagonal entries, we can represent the equation in the form
\begin{equation}
\Psi_{k}(t)
=
\begin{cases}
\mu\nu e^{\mu t}\int\limits_{0}^{\infty}\int\limits_{0}^{\infty}\Psi_{k-1}(u-v)e^{-\mu u-\nu v}dudv,			& \text{if $ t\leq0 $}, \\
\mu\nu e^{-\nu t}\int\limits_{0}^{\infty}\int\limits_{0}^{\infty}\Psi_{k-1}(u-v)e^{-\mu u-\nu v}dudv+1-e^{-\nu t},	& \text{if $ t>0 $}.
\end{cases}\label{E-Psi_k}
\end{equation}

Let us consider the function $ \Psi_{1} $. With the condition $ Y(0)=0 $, we have
$$
\Psi_{1}(t)
=
\mathbb{P}\{\beta_{1}-\alpha_{1}<t\}
=
\begin{cases}
	\frac{\nu}{\mu+\nu}e^{\mu t},			& \text{if $ t\leq0 $}, \\
	1-\frac{\mu}{\mu+\nu}e^{-\nu t},	& \text{if $ t>0 $}.
\end{cases}
$$

Now we take \eqref{E-Psi_k} to get $ \Psi_{2} $. First, we evaluate the integral
$$
\mu\nu\int_{0}^{\infty}\int_{0}^{\infty}\Psi_{1}(u-v)e^{-\mu u-\nu v}dudv
=
\frac{\nu^{2}(3\mu+\nu)}{(\mu+\nu)^{3}}.
$$

Finally, we have
$$
\Psi_{2}(t)
=
\begin{cases}
\frac{\nu^{2}(3\mu+\nu)}{(\mu+\nu)^{3}} e^{\mu t},		& \text{if $ t\leq0 $}, \\
1-\frac{\mu^{2}(\mu+3\nu)}{(\mu+\nu)^{3}} e^{-\nu t},	& \text{if $ t>0 $}.
\end{cases}
$$

\section{Limiting Distribution Function $ \Psi $}

Let us verify the convergence of the sequence $ \Psi_{k} $. Consider the difference
$$
\Psi_{2}(t)-\Psi_{1}(t)
=
\begin{cases}
(\nu-\mu)\frac{\nu\mu}{(\mu+\nu)^{3}} e^{\mu t},		& \text{if $ t\leq0 $}, \\
(\nu-\mu)\frac{\nu\mu}{(\mu+\nu)^{3}} e^{-\nu t},		& \text{if $ t>0 $}.
\end{cases}
$$

For all $ t $, we have either the inequality $ \Psi_{2}(t)-\Psi_{1}(t)\leq0 $ when $ \nu<\mu $, or $ \Psi_{2}(t)-\Psi_{1}(t)\geq0 $ when $ \nu>\mu $. If $ \nu=\mu $, then the functions $ \Psi_{1} $ and $ \Psi_{2} $ coincide.

Assume that $ \nu<\mu $. If $ t\leq0 $, then we have
$$
\Psi_{k+1}(t)-\Psi_{k}(t)
=
\int_{0}^{\infty}\int_{0}^{\infty}(\Psi_{k}(u-v)-\Psi_{k-1}(u-v))f(u-t)g(v)dudv.
$$

It is easy to see by induction that for all $ k\geq1 $, the following inequality holds:
$$
\Psi_{k+1}(t)-\Psi_{k}(t)\leq0.
$$

In the same way, we verify that the inequality is also valid if $ t>0 $.

Considering that $ \Psi_{k}(t)\geq0 $ for all $ t $, we conclude that the sequence $ \Psi_{k}(t) $ converges to some function $ \Psi(t) $ as $ k\to\infty $.

Similar reasoning shows that the sequence of the functions also converges when $ \nu>\mu $. Clearly, if $ \mu=\nu $, then the sequence $ \Psi_{k} $ does not depend on $ k $ at all.

Substitution of the limiting function $ \Psi(t) $ into \eqref{E-Psi_k} gives the equation
\begin{equation}
\Psi(t)
=
\begin{cases}
\mu\nu e^{\mu t}\int\limits_{0}^{\infty}\int\limits_{0}^{\infty}\Psi(u-v)e^{-\mu u-\nu v}dudv,			& \text{if $ t\leq0 $}, \\
\mu\nu e^{-\nu t}\int\limits_{0}^{\infty}\int\limits_{0}^{\infty}\Psi(u-v) e^{-\mu u-\nu v}dudv+1-e^{-\nu t},	& \text{if $ t>0 $}.
\end{cases}\label{E-Psi}
\end{equation}

Taking into account the form of the equation, we find solutions such that
$$
\Psi(t)
=
\begin{cases}
c_{1}e^{\mu t},			& \text{if $ t\leq0 $}, \\
1-c_{2}e^{-\nu t},	& \text{if $ t>0 $};
\end{cases}
$$
where $ c_{1},c_{2} $ are some coefficients.

Let us substitute the above solution into \eqref{E-Psi}. First, we evaluate the integral
$$
\mu\nu\int_{0}^{\infty}\int_{0}^{\infty}\Psi(u-v)e^{-\mu u-\nu v}dudv
=
\frac{\mu\nu}{(\mu+\nu)^{2}}c_{1}
-
\frac{\mu\nu}{(\mu+\nu)^{2}}c_{2}
+
\frac{\nu}{\mu+\nu}.
$$

Now we arrive at the system of algebraic equations with respect to $ c_{1} $ and $ c_{2} $
\begin{align*}
c_{1}
&=
\frac{\mu\nu}{(\mu+\nu)^{2}}c_{1}
-
\frac{\mu\nu}{(\mu+\nu)^{2}}c_{2}
+
\frac{\nu}{\mu+\nu}, \\
c_{2}
&=
-\frac{\mu\nu}{(\mu+\nu)^{2}}c_{1}
+
\frac{\mu\nu}{(\mu+\nu)^{2}}c_{2}
+
\frac{\mu}{\mu+\nu},
\end{align*}
which has the solution
$$
c_{1}
=
\frac{\nu^{2}}{\mu^{2}+\nu^{2}},
\qquad
c_{2}
=
\frac{\mu^{2}}{\mu^{2}+\nu^{2}}.
$$

The limiting function takes the form
$$
\Psi(t)
=
\begin{cases}
	\frac{\nu^{2}}{\mu^{2}+\nu^{2}}e^{\mu t},			& \text{if $ t\leq0 $}, \\
	1-\frac{\mu^{2}}{\mu^{2}+\nu^{2}}e^{-\nu t},	& \text{if $ t>0 $}.
\end{cases}
$$

It is clear that $ \Psi $ presents distribution function of some random variable.

\section{Evaluation of the Mean Growth Rate}

Let us evaluate the mean growth rate $ \lambda $ on the basis of \eqref{E-lambda1}. By applying \eqref{E-Phi}, we get the distribution function
$$
\Phi(t)
=
1
-
\frac{1}{\mu^{2}+\nu^{2}}
\left(
\frac{\mu^{2}+\mu\nu+\nu^{2}}{\mu+\nu}(\nu e^{-\mu t}+\mu e^{-\nu t})
+
\mu\nu e^{-(\mu+\nu)t}
\right).
$$

Its related density function takes the form
$$
\phi(t)
=
\frac{\mu\nu}{\mu^{2}+\nu^{2}}
\left(
\frac{\mu^{2}+\mu\nu+\nu^{2}}{\mu+\nu}(e^{-\mu t}+e^{-\nu t})
-
(\mu+\nu)e^{-(\mu+\nu)t}
\right).
$$

Finally, integrating gives
$$
\lambda
=
\int_{0}^{\infty}t\phi(t)dt
=
\frac{\mu^{4}+\mu^{3}\nu+\mu^{2}\nu^{2}+\mu\nu^{3}+\nu^{4}}{\mu\nu(\mu+\nu)(\mu^{2}+\nu^{2})}.
$$

Note that if $ \nu=\mu $, then we have the result $ \lambda=5/(4\mu) $, which is quite consistent with that obtained in \cite{Jeanmarie1994Analytical}.

\bibliographystyle{utphys}

\bibliography{Growth_rate_of_the_state_vector_in_a_generalized_linear_stochastic_system_with_symmetric_matrix}

\end{document}